\documentclass{article}

\usepackage{arxiv}

\usepackage[utf8]{inputenc} 
\usepackage[T1]{fontenc}    
\usepackage{hyperref}       
\usepackage{url}            
\usepackage{booktabs}       
\usepackage{amsfonts}       
\usepackage{nicefrac}       
\usepackage{microtype}      
\usepackage{lipsum}
\usepackage{graphicx}
\usepackage{amsthm,amsmath,amssymb,color,multirow, inputenc}
\usepackage{wrapfig, float,array}
\newtheorem{thm}{Theorem}[section]

\newtheorem{lem}[thm]{Lemma}

\theoremstyle{definition}

\newtheorem{rem}[thm]{Remark}

\numberwithin{equation}{section}

\graphicspath{ {./images/} }

\title{Non-local boundary value problem for a mixed-type
	equation involving the bi-ordinal Hilfer
	fractional differential operators}

\author{
 Erkinjon Karimov \\
  V.I.Romanovskiy Institute of Mathematics\\
  Uzbekistan Academy of Sciences\\
  Tashkent 100174 \\
  \texttt{erkinjon@gmail.com} \\
   \And
 Bakhodirjon Toshtemirov \\
  V.I.Romanovskiy Institute of Mathematics\\
  Uzbekistan Academy of Sciences\\
  Tashkent 100174 \\
  \texttt{toshtemirovbh@gmail.com} \\
}

\begin{document}
\maketitle
\begin{abstract}
In this paper, we consider a non-local boundary-value problem for a mixed-type equation involving the bi-ordinal Hilfer fractional derivative in rec\-tangular domain.
The main target of this work is to analyze the uniqueness  and the existence of the solution of the considered problem by means of eigenfunctions. Moreover, we construct the solution of the  ordinary fractional  differential equation with the right-sided bi-ordinal Hilfer derivative by the method of reduction to the Volterra integral equation. Then, we present sufficient conditions for given data in order to show the existence of the solution.
 
\end{abstract}


\section{Introduction}
Our project is a competition on Kaggle (Predict Future Sales). We are provided with daily historical sales data (including each products’ sale date, block ,shop price and amount). And we will use it to forecast the total amount of each product sold next month. Because of the list of shops and products slightly changes every month. We need to create a robust model that can handle such situations.

Fractional differential equations have become an important target for investigations because of their properties that are very useful to describe  memory phenomena in control theory, viscoelasti\-city \cite{L1}, anomalous transport and anomalous diffusion\cite{L2}, modeling physical and biological processes \cite{L3}.

Problems for mixed-type equations involving fractional differential operators are often  used to connect two different processes which have different characters on a certain time.
It is known that the theory of boundary-value problems for fractional differential equations is one of the vital fields of the fractional calculus in terms of modeling real-life problems \cite{L4}-\cite{L5}.   We note a few papers devoted to studying boundary-value problems for time-fractional mixed-type partial differential equations \cite{L6}-\cite{L8}.  Of course, one can show some papers regarding to study boundary-value problems for the mixed-type equation involving other fractional differential operators \cite{L9}, \cite{L10} but frankly speaking, the main objects of these investigations were the Riemann-Liouville or the Caputo fractional differential operators.

In 2000, R. Hilfer introduced the generalized Riemann-Liouville derivative (called Hilfer derivative by many authors later) which can be interpolated between Riemann-Liouville and Caputo derivatives \cite{L11}. Several researches have been carried out since that time involving the Hilfer fractional differential operator  \cite{L12}, \cite{L13}, \cite{L13a}  and a few applications  of this operator were investigated by numerous mathematicians. Recently,
V.M. Bulavatsky also generalized the Hilfer derivative in terms of its interpolation concept between the Caputo and the Riemann-Liouvile derivatives of different orders \cite{L14}.

Here we would like to remind another earlier definition of the fractional differential operator that is the  generalization of the Riemann-Liouville, Caputo and Hilfer fractional derivatives introduced by M. M. Dzherbashian and A. B. Nersesian  in 1968 in the following form \cite{L15}
$$
D_{0x}^{\sigma_n}g(x)=I_{0x}^{1-\gamma_n}D_{0x}^{\gamma_{n-1}}...D_{0x}^{\gamma_{1}}D_{0x}^{\gamma_{0}}, ~ ~ n\in\mathbb{N}, ~ x>0
$$
where $I_{0x}^{\alpha}$ and $D_{0x}^{\alpha}$ are the Riemann-Liouville fractional
integral and the Riemann-Liouville fractional derivative of order $\alpha$ respectively,  $\sigma_n\in(0, n]$ which is defined by
$$
\sigma_n=\sum\limits_{j=0}^{n}\gamma_j-1>0, ~ \gamma_j\in(0, 1].
$$
It is called the Dzherbashian-Nersesian differential operator which is less known to the mathematicians than other operators and last year that article was translated and published in FCAA \cite{L16}. The study of problems with this operator was done by many scientists (see \cite{L17}, \cite{L18}). 

In \cite{L19}, the mixed type equation considered where the parabolic part involved the bi-ordinal Hilfer derivative and hyperbolic part with wave equation.

In the present work, we consider the following mixed-type equation with fractional diffusion and fractional wave equation in both parts of the domain involving bi-ordinal Hilfer derivative:

\begin{equation}\label{eq1}
	f(x, t)=\left\{\begin{array}{l} L_1u\equiv{D}_{0+}^{(\alpha_1, \beta_1)\mu_1}u(x, t)-u_{xx}(x, t),\,\,t>0\\ L_2u\equiv{D}_{0-}^{(\alpha_2, \beta_2)\mu_2}u(x, t)-u_{xx}(x, t),\,\,t<0\end{array}\right.
\end{equation}
in mixed domain $\Omega=\Omega_1\cup\Omega_2\cup AB$. Where $f(x, t)$ is a given function,\\  $i-1<\alpha_i, \beta_i<{i}$, \,  $0\leq\mu_i\leq1$, $i=\overline{1,2}$, \,  $\Omega_1=\left\{(x,t):\, 0<x<l,\,0<t<T\right\}$, ~  $\Omega_2=\left\{(x,t):\,0<x<l,\,-T<t<0\right\}$, ~ $T>0$, ~ $AB=\left\{(x,t):\,0<x<l,\,t=0\right\}$,
\begin{equation}\label{e2}
	D_{0\pm}^{(\alpha_i, \beta_i)\mu_i}g(t)=I_{0\pm}^{\mu_i(i-\alpha_i)}\left(\pm\frac{d}{dt}\right)^iI_{0\pm}^{(1-\mu_i)(i-\beta_i)}g(t)
\end{equation}
is the bi-ordinal Hilfer fractional differential operator (FDO) of orders $\alpha_i, \beta_i $  and of type $\mu_i$, where
\[
I_{0+}^{\alpha}g(t)=\frac{1}{\Gamma(\alpha)}\int\limits_0^t\frac{g(z)dz}{(t-z)^{1-\alpha}},
\]
\[
I_{0-}^{\alpha}g(t)=\frac{1}{\Gamma(\alpha)}\int\limits_t^0\frac{g(z)dz}{(z-t)^{1-\alpha}}
\]
are the left-sided and right-sided Riemann-Liouville fractional integral of order $\alpha\, (\Re (\alpha)>0)$ \cite{L4}, where $\Gamma(\alpha)$ is Euler's gamma-function.

We note that the bi-ordinal Hilfer fractional derivative can be reduced from the Dzherbashian-Nersesian differential operator as a particular case, i.e. for $n=1$
$$
D_{0+}^{\sigma_1}g(x)=I_{0+}^{1-\gamma_1}D_{0+}^{\gamma_0}g(x).
$$

Nonlocal BVP for Eq.(\ref{eq1}) in $\Omega$ can be formulated as follows:

\textbf{Problem B.} Find a solution $u(x,t)$ of equation (\ref{eq1}) which is subject to the following regularity conditions
$$
t^{1-\gamma_1}u(x, t), ~ t^{1-\gamma_1}D_{0+}^{(\alpha_1, \beta_1)\mu_1}u(x, t)\in C(\overline{\Omega}_1), \, \,
(-t)^{2-\gamma_2}u(x, t)\in C(\overline{\Omega}_2),
$$
$$
~ (-t)^{2-\gamma_2}D_{0-}^{(\alpha_2, \beta_2)\mu_2}u(x, t)\in C(\overline{\Omega}_2),
~ u_{xx}\in C(\Omega_1\cup\Omega_2),
$$
submitted to the boundary conditions
\begin{equation}\label{eq2}
	u(0,t)=0,\,\,u(l,t)=0, ~  t\in[-T,0)\cap(0, T]
\end{equation}
and non-local condition
\begin{equation}\label{eq3}
	u(x,-T)=u(x, T)+\psi(x),\,\,\,0\leq x\leq l,
\end{equation}
and also it satisfies the conjugation conditions on $AB$
\begin{equation}\label{eq4}
	\lim\limits_{t\rightarrow +0}I_{0+}^{1-\gamma_1}u(x, t)= \lim\limits_{t\rightarrow -0}I_{0-}^{2-\gamma_2}u(x, t),\,\,0\leq x\leq l,
\end{equation}
\begin{equation}\label{eq5}
	\begin{array}{l}
		\lim\limits_{t\rightarrow +0} t^{1-\delta_1}\left(\frac{\partial}{\partial{t}}I_{0+}^{1-\gamma_1}u(x,t)\right)=\lim\limits_{t\rightarrow -0}\frac{\partial}{\partial{t}}I_{0-}^{1-\gamma_2}u(x, t) \,\,0< x< l.
	\end{array}
\end{equation}
Where $\gamma_i=\beta_i+\mu_i(i-\beta_i)$ \, $\delta_i=\beta_i+\mu_i(\alpha_i-\beta_i)$, \,  $(i=\overline{1,2})$, $\psi(x)$ is a given function such that $\psi(0)=\psi(l)=0$.

The rest of the article is organized as follows. In Section 2 we recall some definitions and solve the Cauchy problem with the right-sided bi-ordinal Hilfer derivative which are needed in the sequel. Main results is presented in Section 3. Furthermore, we analyze the uniqueness of the solution by using eigenfunctions and  we present sufficient conditions for given data in order to show the existence of the result.

\section{Preliminaries}

The two parameter Mittag-Leffler function \cite{L4} is defined as

$$
E_{\alpha, \beta}(z)=\sum_{k=0}^{\infty}\frac{z^k}{\Gamma(\alpha{k}+\beta)}, \, \, \, \alpha>0, \, \beta\in{R}.
$$

\begin{lem}\label{lm1}
 (see \cite{L4}) Let $\alpha<2,$ \, $\beta$ is an arbitrary real number, and ${\pi\alpha}/{2}<\mu<min\{\pi, \pi\alpha\}$, such that $\mu\leq|argz|\leq\pi, \, \, |z|\geq0$, $M$ is a positive constant. Then, the following estimate hold
$$
|E_{\alpha, \beta}(z)|\leq\frac{M}{1+|z|}.
$$
\end{lem}
If $\alpha>{0}$,  $\nu>0$ and $\beta>0$, then the following results hold \cite{L29}:
\begin{equation}\label{ml1}
	E_{\alpha, \beta}(z)=\frac{1}{\Gamma(\beta)}+zE_{\alpha, \alpha+\beta}(z),
\end{equation}
\begin{equation}\label{ml2}
	\frac{1}{\Gamma(\nu)}\int\limits_{t}^{0}(z-t)^{\nu-1}E_{\alpha, \beta}(\lambda (-z)^{\alpha})(-z)^{\beta-1}dz=(-t)^{\beta+\nu-1}E_{\alpha, \beta+\nu}(\lambda (-t)^{\alpha}).
\end{equation}

We recall the following lemma which is analogy of the lemma presented by A. Pskhu \cite{L20}.

\begin{lem}\label{lm2} Let $\alpha>0$ and $\pi\geq|argz|>\frac{\pi\alpha}{2}+\varepsilon, \, \varepsilon>0$, then the following relations are valid for $|z| \to +\infty$
$$
\lim\limits_{|z|\rightarrow +\infty}E_{\alpha, \beta}(z)=0,
$$
$$
\lim\limits_{|z|\rightarrow +\infty}zE_{\alpha, \beta}(z)=-\frac{1}{\Gamma(\beta-\alpha)}.
$$
\end{lem}

The Riemann-Liouville fractional derivatives left-sided  $D^{\alpha}_{0+}f(t)$ and the right-sided $D^{\alpha}_{0-}f(t)$ of order $\alpha$  $(n-1<\alpha\leq{n})$  are defined as \cite{L5}

$$
D^{\alpha}_{0+}g(t)=\left(\frac{d}{dt}\right)^nI^{n-\alpha}_{0+}g(t), \, \, \, \,   n-1\leq\alpha<n,
$$
$$
D^{\alpha}_{0-}g(t)=(-1)^n\left(\frac{d}{dt}\right)^nI^{n-\alpha}_{0-}g(t), \, \, \, \,   n-1\leq\alpha<n,
$$

respectively.

We rewrite some basic properties of the right-sided Riemann-Liouville integral and differential operators(see \cite{L5}):

1) If $\alpha>0$, $\beta>0$ and $t\in[-T, 0]$, $f(t)\in{L}_p(-T, 0), (1\leq{p}<\infty)$
$$
I_{0-}^{\alpha}I_{0-}^{\beta}g(t)=I_{0-}^{\alpha+\beta}g(t)
$$

2) If $g(t)\in{L_1}(-T, 0)$ and $I^{n-\alpha}_{0-}g(t)\in{AC}^{n}[-T, 0]$ then the inequality
\begin{equation}\label{ee7}
	I_{0-}^{\alpha}D_{0-}^{\alpha}g(t)=g(t)-\sum_{j=1}^{n}\frac{(-1)^{n-j}(-t)^{\alpha-j}}{\Gamma(\alpha-j+1)}\left[\mathop {\lim }\limits_{t \to  0-}\left(\frac{d}{dt}\right)^{n-j}I_{0-}^{n-\alpha}g(t)\right]
\end{equation}

hods almost everywhere on $[-T, 0]$, where
$$
AC^{n}[-T, 0]=\left\{g:[-T, 0]\to \mathbb{C} ~ \text{and} ~ \left(\frac{d}{dx}\right)^{n-1}g(x)\in AC[-T, 0], \right\},
$$
$$
g(x)\in AC[-T, 0]\Leftrightarrow g(x)=c+\int_{-T}^{x}\varphi(t)dt, (\varphi(t)\in L(-T, 0))
$$
The left-sided $^CD_{0+}^{\alpha}f(t)$ and right-sided $^CD_{0-}^{\alpha}f(t)$ Caputo derivatives of order $\alpha \, \, (n-1<\alpha\leq{n})$ are defined by (see \cite{L5})
$$
^CD_{0+}^{\alpha}g(t)=D_{0+}^{\alpha}\left[g(t)-\sum_{k=0}^{n-1}\frac{g^{(k)}(0+)}{k!}t^k\right]
$$
and
$$
^CD_{0-}^{\alpha}g(t)=D_{0-}^{\alpha}\left[f(t)-\sum_{k=0}^{n-1}\frac{g^{(k)}(0-)}{k!}(-t)^k\right]
$$
respectively.
\begin{rem}
We have the following comments:

1) The bi-ordinal Hilfer derivative $	D^{(\alpha, \beta) \mu}_{0\pm}g(t)$ can be written as
$$
D_{0\pm}^{(\alpha, \beta)\mu}g(t)=I^{\mu(n-\alpha)}_{0\pm}\left(\pm\frac{d}{dt}\right)^nI^{(1-\mu)(n-\beta)}_{0\pm}g(t)=
$$
\begin{equation}\label{ee8}
	=I_{0\pm}^{\mu(n-\alpha)}\left(\pm\frac{d}{dt}\right)^nI_{0\pm}^{n-\gamma}g(t)=I^{\mu(n-\alpha)}_{0\pm}D_{0\pm}^{\gamma}g(t)=I_{0\pm}^{\gamma-\delta}D_{0\pm}^{\gamma}g(t),
\end{equation}
for $t\in[-T, T]$, where $\gamma=\beta+\mu(n-\beta)$ and $\delta=\beta+\mu(\alpha-\beta)$.

2) In general, (\ref{e2}) is also preserved in terms of its interpolation concept. Specifically, when $\mu=0$, (\ref{e2}) denotes Riemann-Liouville fractional derivative operator of $\beta$ order and for $\mu=1$, the bi-ordinal Hilfer fractional derivative (\ref{e2}) expresses the Caputo fractional derivative of order $\alpha$ i.e.
$$
D_{a+}^{(\alpha, \beta)\mu}g(t)=\left\{ \begin{gathered}
	D_{a+}^{\beta}g(t), \, \, \, \, \mu=0,\hfill \cr
	^CD_{a+}^{\alpha}g(t), \, \, \, \, \, \mu=1. \hfill \cr
\end{gathered}  \right.
$$
\end{rem}
\begin{lem}\label{lm3}
 If $f(t)\in L_1(-T, 0)$, $T>0$, $\alpha>0$, $\lambda\in\mathbb{C}$, then
$$
y(t)-\frac{\lambda}{\Gamma(\alpha)} \int\limits_{t}^{0}\left(s-t\right)^{\alpha-1}E_{\alpha, \alpha}\left[\lambda(s-t)^{\alpha}\right]y(s)ds=g(t),
$$
integral equation has the following unique solution
$$
y(t)=g(t)+\lambda \int\limits_{t}^{0}(s-t)^{\alpha-1}E_{\alpha, \alpha}[\lambda(s-t)^{\alpha}]g(s)ds.
$$
\end{lem}

The solution of the integral equation with the left-sided Riemann-Liouville integral operator second kind   was found in \cite{L23} by means of the Laplace transform.
For ascertaining the result of Lemma \ref{lm3} one can easily check by substituting the result into the equation.
\subsection*{Cauchy problem for the equation with the right-sided bi-ordinal Hilfer fractional derivative}

Let us consider the problem about finding a solution of  the following ordinary differential equation involving the right-sided bi-ordinal Hilfer fractional derivative with the initial conditions:
\begin{equation}\label{eq8}
	\left\{\begin{array}{l} D_{0-}^{(\alpha, \beta)\mu}u(t)=\lambda u(t)+g(t),\\ \lim\limits_{t\to 0-}I_{0-}^{2-\gamma}u(t)=\xi_0,\\ \lim\limits_{t\to 0-}\frac{d}{dt}I_{0-}^{2-\gamma}u(t)=\xi_1,
	\end{array}\right.
\end{equation}
where $1<\alpha, \beta\leq2$, $\gamma=\beta+\mu(2-\beta)$, $\xi_0, \xi_1\in \mathbb{R},$ $g(t)$ is the given function.

\begin{lem}\label{lm4}
Let $g(t)\in C[-T, 0], ~ g'(t)\in L_{1}(-T, 0)$. Then the solution of the problem (\ref{eq8}) as follows:
$$
u(t)=\xi_0(-t)^{\gamma-2}E_{\delta, \gamma-1}[\lambda(-t)^{\delta}]-\xi_1(-t)^{\gamma-1}E_{\delta, \gamma}[\lambda(-t)^{\delta}]+
$$
\begin{equation}\label{e9}
	+\int\limits_{t}^{0}(z-t)^{\delta-1}E_{\delta, \delta}\left[\lambda(z-t)^{\delta}\right]g(z)dz,
\end{equation}
where $\delta=\beta+\mu(\alpha-\beta)$,  $\gamma=\beta+\mu(2-\beta)$.
\end{lem}

\noindent{\bf Proof.} First, we rewrite the equation according to (\ref{ee8}) as follows:
$$
I_{0-}^{\gamma-\delta}D_{0-}^{\gamma}u(t)=\lambda u(t)+g(t).
$$
By applying $I_{0-}^{\delta}$ operator to this equation we get
$$
I_{0-}^{\gamma}D_{0-}^{\gamma}u(t)=\lambda I_{0-}^{\delta}u(t)+I_{0-}^{\delta}g(t)
$$
and using (\ref{ee7}) properties of the Riemann-Liouville integral and differential operators yield
$$
u(t)=\lambda I_{0-}^{\delta}u(t)+I_{0-}^{\delta}g(t)+\frac{(-t)^{\gamma-2}\xi_0}{\Gamma(\gamma-1)}-\frac{(-t)^{\gamma-1}\xi_1}{\Gamma(\gamma)},
$$
This integral equation has the following solution according to Lemma \ref{lm3}
$$
u(t)=g^*(t)+\lambda \int\limits_{t}^{0}(s-t)^{\delta-1}E_{\delta, \delta}[\lambda(s-t)^{\delta}]g^*(s)ds=L_1(t)+L_2(t),
$$
where $g^*(t)=I_{0-}^{\delta}g(t)+\frac{(-t)^{\gamma-2}\xi_0}{\Gamma(\gamma-1)}-\frac{(-t)^{\gamma-1}\xi_1}{\Gamma(\gamma)}$.
$$
L_1(t)=\frac{(-t)^{\gamma-2}\xi_0}{\Gamma(\gamma-1)}-\frac{(-t)^{\gamma-1}\xi_1}{\Gamma(\gamma)}+
$$
$$
+\lambda \int\limits_{t}^{0}(s-t)^{\delta-1}E_{\delta, \delta}[\lambda(s-t)^{\delta}]\Big(\frac{(-s)^{\gamma-2}\xi_0}{\Gamma(\gamma-1)}-\frac{(-s)^{\gamma-1}\xi_1}{\Gamma(\gamma)}\Big)ds,
$$
$$
L_2(t)=I_{0-}^{\delta}g(t)+\lambda \int\limits_{t}^{0}(s-t)^{\delta-1}E_{\delta, \delta}[\lambda(s-t)^{\delta}]I_{0-}^{\delta}g(s)ds.
$$

If we use $z=t-s$ substitution to the integral in $L_1(t)$  and using (\ref{ml1}), (\ref{ml2}) properties we can easily obtain the following result
\begin{equation}\label{ee13}
	L_1(t)=\xi_0(-t)^{\gamma-2}E_{\delta, \gamma-1}[\lambda(-t)^{\delta}]-\xi_1(-t)^{\gamma-1}E_{\delta, \gamma}[\lambda(-t)^{\delta}].
\end{equation}
Now we consider second integral on $L_2(t)$
$$
\int\limits_{t}^{0}(s-t)^{\delta-1}E_{\delta, \delta}[\lambda(s-t)^{\delta}]I_{0-}^{\delta}f(s)ds=
$$
$$
=\frac{1}{\Gamma(\delta)}\int\limits_{t}^{0}(s-t)^{\delta-1}E_{\delta, \delta}[\lambda(s-t)^{\delta}]ds\int\limits_{s}^{0}(z-s)^{\delta-1}f(z)dz=
$$
$$
=\frac{1}{\Gamma(\delta)}\int\limits_{t}^{0}f(z)dz\int\limits_{t}^{z}(z-s)^{\delta-1}(s-t)^{\delta-1}E_{\delta, \delta}[\lambda(s-t)^{\delta}]ds.
$$
By means of formula (\ref{ml2}) we simplify the integrant as follows
$$
\int\limits_{t}^{z}(z-s)^{\delta-1}(s-t)^{\delta-1}E_{\delta, \delta}[\lambda(s-t)^{\delta}]ds=\Gamma(\delta)(z-t)^{2\delta-1}E_{\delta, 2 \delta}(\lambda(z-t)^{\delta}).
$$
To clarify further, we use (\ref{ml1}), then the form of $L_2(t)$ will be as follows
\begin{equation}\label{ee14}
	L_2(t)=\int\limits_{t}^{0}(z-t)^{\delta-1}E_{\delta, \delta}\left[\lambda(z-t)^{\delta}\right]f(z)dz.
\end{equation}

Finally, we can obtain the solution presented in Lemma \ref{lm4} from (\ref{ee13}) and (\ref{ee14}). The similar lemma to Lemma \ref{lm4} was also studied in \cite{L13a} for equation with the left-sided Hilfer differential operator. The proof of Lemma \ref{lm4} is completed.
\section{Main Results}

We intend to prove the uniqueness and existence of solution to the problem for the mixed-type equation (\ref{eq1}) along  with the conditions (\ref{eq2})-(\ref{eq5}).

\subsection{Investigation  of the main problem}

 First let us introduce the following new  notations
 \begin{equation}\label{eq6}
 	\tau(x)=\lim\limits_{t\rightarrow+0}I_{0+}^{1-\gamma_1}u(x,t),\,\,0\leq x\leq l,
 \end{equation}
 \begin{equation}\label{e8}
 	\varphi(x)=\lim\limits_{t\rightarrow -0}I_{0-}^{2-\gamma_2}u(x,t),\,\,0\leq x\leq l,
 \end{equation}
 \begin{equation}\label{e9}
 	\nu(x)=\lim\limits_{t\rightarrow -0}\frac{\partial}{\partial{t}}I_{0-}^{2-\gamma_2}u(x,t),\,\,0<x<l.
 \end{equation}

 For solving the problem we use the method of separation of variables for homogeneous equation corresponding (\ref{eq1}) along with the conditions (\ref{eq2}) and we obtain the spectral problem which is its eigenvalues and eigenfunctions are in the following forms
 \begin{equation}\label{sp}
 	\lambda_n=\left(\frac{n\pi}{l}\right)^2, \, \, \, X_n(x)=\sin(\sqrt{\lambda_n}x).
 \end{equation}

 The system of $X_n(x)$ in the form (\ref{sp} is the orthogonal basis in $L_2(0, l)$, for that reason we can represent the solution $u(x, t)$ and the given function $f(x, t)$ in the form of series expansions as follows
 \begin{equation}\label{u}
 	u(x, t)=\sum\limits_{n=1}^{\infty}u_n(t)\sin(\sqrt{\lambda_n}x)
 \end{equation}
 and
 \begin{equation}\label{f}
 	f(x, t)=\sum\limits_{n=1}^{\infty}f_n(t)\sin(\sqrt{\lambda_n}x),
 \end{equation}
 where

 \begin{equation}\label{fn}
 	f_n(t)=\frac{2}{l}\left\{\begin{array}{l} \int\limits_{0}^{t}f(x, t)\sin(\sqrt{\lambda_n}x),  \,\,t>0,\\ \int\limits_{t}^{0}f(x, t)\sin(\sqrt{\lambda_n}x), \,\,t<0.\end{array}\right.
 \end{equation}

 Substituting (\ref{u}) and (\ref{f}) into the equation (\ref{eq1}) along with the conditions (\ref{eq6}) and (\ref{e8}), (\ref{e9}) we obtain the problem for the ordinary fractional differential equations in $\Omega_1$ and $\Omega_2$ respectively.
 
 We note works of A. Pskhu \cite{Lp24}, \cite{Lp25}, where main boundary value problems for diffusion-wave equation with the Riemann-Liouville fractional derivative were investigated by the method of Green's functions. 

 The ordinary fractional differential equation with respect to $t$ corresponding Eq.(\ref{eq1}) has been studied in \cite{L14} for $t>0$. Hence, we can write the solution of the Eq.(\ref{eq1}) in $\Omega_1$ which satisfies conditions (\ref{eq2}), (\ref{eq6}) as follows:
 \begin{equation*}\label{eq7}
 	u(x,t)=\sum\limits_{n=1}^{+\infty}\tau_n t^{\gamma_1-1}E_{\delta_1, \gamma_1}(-\lambda_nt^{\delta_1})\sin(\sqrt{\lambda_n}x)+
 	\end{equation*}
 \begin{equation}
 	+\int\limits_{0}^{t}(t-s)^{\delta_1-1}E_{\delta_1, \delta_1}\left[-\lambda_n(t-s)^{\delta_1}\right]f_n(s)ds\sin(\sqrt{\lambda_n}x),
 \end{equation}
 where $\lambda_n=\left(\frac{n\pi}{l}\right)^2$.

 Using representations (\ref{eq7}), we evaluate $t^{1-\delta_1}\left(I_{0+}^{1-\gamma_1}u(x,t)\right)_t$:
 \[
 t^{1-\delta_1}\left(I_{0+}^{1-\gamma_1}u(x,t)\right)_t=t^{1-\delta_1}\sum_{n=1}^{\infty}\frac{d}{dt}\tau_nE_{\delta_1, 1}(-\lambda_nt^{\delta_1})\sin(\sqrt{\lambda_n}x)+
 \]
 \[
 +t^{1-\delta_1}\sum_{n=1}^{\infty}\frac{d}{dt}\int\limits_0^t(t-s)^{2\delta_1-1}E_{\delta_1, \delta_1-\gamma_1+1}[-\lambda_n(t-s)^{\delta}]f_n(s)ds\sin(\sqrt{\lambda_n}x)=
 \]
 \[
 =t^{1-\delta_1}\sum_{n=1}^{\infty}\left(-\lambda_n\tau_nt^{\delta_1-1}E_{\delta_1, 1}(-\lambda_nt^{\delta_1})-f(0)t^{2\delta_1-1}E_{\delta_1, \delta_1-\gamma+1}(-\lambda t^{\delta_1})\right)\sin(\sqrt{\lambda_n}x)-
 \]
 \[
 -t^{1-\delta_1}\sum_{n=1}^{\infty}\left(\int\limits_0^t(t-s)^{2\delta_1-1}E_{\delta_1, \delta_1-\gamma_1+1}[-\lambda_n(t-s)^{\delta}]f'_n(s)ds\right)\sin(\sqrt{\lambda_n}x).
 \]

 According to the above evaluation, we can calculate the limit
 \begin{equation}\label{eq9}
 	\lim\limits_{t\rightarrow+0}t^{1-\delta_1}\left(I_{0+}^{1-\gamma}u(x,t)\right)_t=\sum\limits_{n=1}^{\infty}(-\lambda_n)\tau_n\sin(\sqrt{\lambda_n}x),\,\,0<x<l.
 \end{equation}

 Considering notations (\ref{eq6}), (\ref{e8}) and conjugation condition (\ref{eq4}), as such from (\ref{e9}), (\ref{eq9}) and conjugation condition (\ref{eq5}), we obtain the following linear equations:
 \begin{equation}\label{e13}
 	\left\{\begin{array}{l} \tau_n=\varphi_n,\\ -\frac{\lambda_n\tau_n}{\Gamma(\delta_1)}=\nu_n,\end{array}\right.
 \end{equation}
where $\tau_n$, $\varphi_n$ and $\nu_n$ are  Fourier coefficients of the unknown functions $\tau(x)$, $\varphi(x)$ and $\nu(x)$ respectively.

 Now we will establish another functional relation which is determine from (\ref{eq3}). For this aim, we need the solution of the  problem intended to solve $L_2u=0$ equation with the conditions (\ref{e8}), (\ref{e9}). After applying method of separation variables, we have spectral problem which its eigenvalues and eigenfunctions as given in (\ref{sp}) and the problem for ordinary differential equation involving the right-sided bi-ordinal Hilfer fractional differential operator as (\ref{eq8}).

 According to Lemma 4 and (\ref{sp}) we can write the solution of $L_2u=0$ satisfying (\ref{eq2}), (\ref{e8}), (\ref{e9}) conditions as follows:
 \[
 u(x, t)=\sum\limits_{n=1}^{+\infty}\varphi_n(-t)^{\gamma_2-2}E_{\delta_2, \gamma_2-1}[-\lambda_n(-t)^{\delta_2}]\sin(\sqrt{\lambda_n}x)-
 \]
 \begin{equation*}
 -\sum\limits_{n=1}^{+\infty}\nu_n(-t)^{\gamma_2-1}E_{\delta_2, \gamma_2}[-\lambda_n(-t)^{\delta_2}]\sin(\sqrt{\lambda_n}x)+	
  \end{equation*}
 \begin{equation}\label{e14}
 	+\sum\limits_{n=1}^{+\infty}\int\limits_{t}^{0}(z-t)^{\delta_2-1}E_{\delta_2, \delta_2}[-\lambda_n(z-t)^{\delta_2}]f_n(z)dz\sin(\sqrt{\lambda_n}x)
 \end{equation}

 Substituting (\ref{eq7}) and (\ref{e14}) into  (\ref{eq3})  we deduce that
 $$
 \psi_n=\varphi_nT^{\gamma_2-2}E_{\delta_2, \gamma_2-1}\left(-\lambda_n{T}^{\delta_2}\right)-\nu_nT^{\gamma_2-1}E_{\delta_2, \gamma_2}(-\lambda_n{T}^{\delta_2})+
 $$
 $$
 +\int\limits_{-T}^{0}(z+T)^{\delta_2-1}E_{\delta_2, \delta_2}(-\lambda(z+T)^{\delta_2})f_n(z)dz	-\tau_nT^{\gamma_1-1}E_{\delta_1, \gamma_1}(-\lambda{T}^{\delta_1})
 $$
  \begin{equation}\label{eq12}
 -\int\limits_{0}^{T}(T-z)^{\delta_1-1}E_{\delta_1, \delta_1}(-\lambda_n(T-z)^{\delta_1})f_n(z)dz.
 \end{equation}

 From the equalities (\ref{e13}) and (\ref{eq12}) we determine $\tau_n, \, \varphi_n, \, \nu_n$ unknowns in the following forms
 \begin{equation}\label{tau}
 	\tau_n=\frac{1}{\Delta_n}(\psi_n+F_n),
 \end{equation}
 \begin{equation}\label{nu}
 	\nu_n=\frac{-\lambda_n}{\Delta_n}(\psi_n+F_n),
 \end{equation}
 \begin{equation}\label{fi}
 	\varphi_n=\frac{1}{\Delta_n}(\psi_n+F_n),
 \end{equation}
 where
 \begin{equation*}
 		\Delta_n=T^{\gamma_2-2}E_{\delta_2, \gamma_2-1}\left(-\lambda_n{T}^{\delta_2}\right)+
 	\end{equation*}
 \begin{equation}\label{e16}
 +\frac{\lambda_nT^{\gamma_2-1}}{\Gamma(\delta_1)}E_{\delta_2, \gamma_2}(-\lambda_n{T}^{\delta_2})-T^{\gamma_1-1}E_{\delta_1, \gamma_1}(-\lambda{T}^{\delta_1}),
 \end{equation}
 $$
 F_n=\int\limits_{0}^{T}(T-z)^{\delta_1-1}E_{\delta_1, \delta_1}(-\lambda_n(T-z)^{\delta_1})f_n(z)dz-
 $$
 $$
 -\int\limits_{-T}^{0}(z+T)^{\delta_2-1}E_{\delta_2, \delta_2}(-\lambda_n(z+T)^{\delta_2})f_n(z)dz.
 $$

 \subsection{A uniqueness of the solution}
 We assume that there exist two  different $u_1(x, t)$ and $u_2(x, t)$ solutions of the main problem. Then it is enough to show that $u(x, t)=u_1(x, t)-u_2(x, t)$ is a trivial solution of the homogeneous problem.

 Let $u(x, t)$ be a solution of the homogeneous problem.

 Let us first consider the following integral
 \begin{equation}\label{eq20}
 	u_n(t)=\int\limits_0^1u(x, t)\sin(\sqrt{\lambda_n}x)dx, \, n=1, 2, 3, ....,
 \end{equation}

 Then we introduce another function based on (\ref{eq20})
 \begin{equation}\label{e18}
 	v^{\varepsilon }_{n}(t)=\int\limits_{\varepsilon}^{1-\varepsilon}u(x, t)\sin(\sqrt{\lambda_n}x)dx, \, n=1, 2, 3, ....,
 \end{equation}

 Applying $D_{0+}^{(\alpha_1, \beta_1)\mu_1}$ and $D_{0-}^{(\alpha_2, \beta_2)\mu_2}$ to (\ref{e18}) and using the equation(\ref{eq1})

 $$
 D_{0+}^{(\alpha_1, \beta_1)\mu_1}v_{\varepsilon}(t)=-2\int\limits_{\varepsilon}^{1-\varepsilon}D_{0+}^{(\alpha_1, \beta_1)\mu_1}u(x, t)\sin(\sqrt{\lambda_n}x)dx=
 $$
 $$
 =-2\int\limits_{\varepsilon}^{1-\varepsilon}u_{xx}(x, t)\sin(\sqrt{\lambda_n}x)dx
 $$

 $$
 D_{0-}^{(\alpha_2, \beta_2)\mu_2}v_{\varepsilon}(t)=-2\int\limits_{\varepsilon}^{1-\varepsilon}D_{0-}^{(\alpha_2, \beta_2)\mu_2}u(x, t)\sin(\sqrt{\lambda_n}x)dx=
 $$
 $$
 =-2\int\limits_{\varepsilon}^{1-\varepsilon}u_{xx}(x, t)\sin(\sqrt{\lambda_n}x)dx
 $$
 and integrating by parts twice the right sides of the equalities on $t\in(0, T)$ and $t\in(-T, 0)$, respectively, and passing to the limit on $\varepsilon \to +0$ yield
 \begin{equation}\label{e19}
 	\left\{\begin{array}{l} D_{0+}^{(\alpha_1, \beta_1)\mu_1}u_{n}(t)+\lambda^2u_n(t)=0,  \, \, \,   t>0, \\ D_{0-}^{(\alpha_2, \beta_2)\mu_2}u_{n}(t)+\lambda^2u_n(t)=0, \, \, \, t<0. \end{array}\right.
 \end{equation}

 Considering conditions (\ref{eq6}), (\ref{e8}), (\ref{e9}) in homogeneous case, (\ref{e19}) has a solution $u_n(t)=0$ if $\Delta_n\neq0$ in (\ref{e16}). Then from (\ref{eq20}) and the completeness of the system $\{X_n(x)\}$ in the space $L_2[0, l]$, $u(x, t)\equiv0$ in $\overline{\Omega}$. This completes the prove of uniqueness of the solution of the main problem.

 \subsection{Existence of a solution}

 First of all, we prove that $\Delta_n\neq0$ for sufficient large $n$. Considering Lemma \ref{lm2} we can show
 $$
 \lim\limits_{n\rightarrow +\infty}\Delta_n=\lim\limits_{\lambda_n\rightarrow +\infty}\Delta_n=\lim\limits_{|z_1|\rightarrow +\infty}\left(T^{\gamma_2-2}E_{\delta_2, \gamma_2-1}(z_1)-\frac{T^{\gamma_2-1-\delta_2}}{\Gamma(\delta_1)}E_{\delta_2, \gamma_2}(z_1)\right)-
 $$

 $$
 -\lim\limits_{|z_2|\rightarrow +\infty}T^{\gamma_1-1}E_{\delta_1, \gamma_1}(z_2)=\frac{T^{\gamma_2-\delta_2-1}}{\Gamma(\delta_1)\Gamma(\gamma_2-\delta_2)}>0.
 $$

 In other words, it confirms that $\Delta_n>0$ for any sufficient large $n$.

 \smallskip

 For showing the existence of the result, we prove the uniform convergence of the series of $u(x, t)$,   $u_{xx}(x, t)$ and $D_{0\pm}^{(\alpha_i, \beta_i)\mu_i}u(x, t)$, $i=\overline{1,2}.$

 First, we get the estimates of the function $u(x, t)$ in $\Omega_1$ by means of Lemma \ref{lm1}:
 $$
 |t^{1-\gamma_1}u(x, t)|\leq\sum\limits_{n=1}^{\infty}|\tau_n|E_{\delta_1, \gamma_1}(-\lambda_nt^{\delta_1})|+
 $$
 $$
 +\sum\limits_{n=1}^{\infty}t^{1-\gamma_1}\int\limits_{0}^{t}|t-s|^{\delta_1-1}|E_{\delta_, \delta_1}[-\lambda_n(t-s)^{\delta_1}]||f_n(s)|ds=
 $$
 $$
 =\sum\limits_{n=1}^{\infty}|\tau_n|E_{\delta_1, \gamma_1}(-\lambda_nt^{\delta_1})|+ |f_n(0)||t|^{\delta_1}|E_{\delta_1, \delta_1+1}(-\lambda_nt^{\delta_1})|+
 $$
 $$
 +\sum\limits_{n=1}^{\infty}|t|^{1-\gamma_1}\int\limits_{0}^{t}|t-s|^{\delta_1}|E_{\delta_, \delta_1+1}[-\lambda_n(t-s)^{\delta_1}]||f'_n(s)|ds=
 $$
 $$
 \leq\sum\limits_{n=1}^{\infty}\left(\frac{|\tau_n||M}{1+\lambda_n|t^{\delta_1}|}
 +\frac{|f_n(0)||t|^{\delta_1+1-\gamma_1}}{1+\lambda_n|t|^{\delta_1}}+T^{1-\gamma_1}\int\limits_{0}^{t}\frac{|t-s|^{\delta_1}M}{1+\lambda_n|t-s|^{\delta_1}}|f'_n(s)|ds\right).
 $$
 If $\tau(x)\in{C}\, [0,l], ~ \tau'(x)\in L_2(0, l)$ and $f(x, t)\in C^{0,1}[0, l]\times [0, T]$,   then the series of $u(x, t)$ is bounded with convergent numerical series. Note that the condition $\tau'(x)\in L_2(0, l)$ is required for $t\to+0$.  From Weierstrass M-test the series of $u(x, t)$ is considered uniformly convergent in $\Omega_1$.

 As such for estimate $u(x, t)$ in $\Omega_2$ after integrating by parts

 \begin{equation*}
 	|(-t)^{2-\gamma_2}u(x, t)|\leq\sum\limits_{n=1}^{+\infty}\left(\varphi_n \Big|E_{\delta_2, \gamma_2-1}[-\lambda_n(-t)^{\delta_2}]\Big|+\nu_n|-t|\Big|E_{\delta_2, \gamma_2}[-\lambda_n(-t)^{\delta_2}]\Big|\right)+
 \end{equation*}
 \begin{equation*}
 	+\sum\limits_{n=1}^{\infty}|-t|^{2-\gamma_2+\delta_2}|f_n(0)||E_{\delta_2, \delta_2+1}(-\lambda_n(-t)^{\delta_2})|+
 \end{equation*}
 \begin{equation*}
 	+\sum\limits_{n=1}^{\infty}|-t|^{2-\gamma_2}\int\limits_{t}^{0}|z-t|^{\delta_2}\Big|E_{\delta_2, \delta_2+1}[-\lambda_n(z-t)^{\delta_2}]\Big| |f'_n(z)|dz=
 \end{equation*}
 \begin{equation*}
 	\leq \sum\limits_{n=1}^{\infty}\left(\frac{\varphi_n M}{1+\lambda_n |-t|^{\delta_2}}+\frac{\nu_n |-t| M}{1+\lambda_n |-t|^{\delta_2}}\right)+
 \end{equation*}
 \begin{equation*}	
 	+\sum\limits_{n=1}^{\infty}\left(\frac{|f_n(0)||-t|^{\delta_2} M}{1+\lambda_n |-t|^{\delta_2}}+T^{2-\gamma_2}\int\limits_{t}^{0}\frac{|z-t|^{\delta_2} M}{1+\lambda_n|z-t|^{\delta_2}} |f'_n(z)|dz\right).
 \end{equation*}

 If $\varphi(x), ~ \nu(x)\in C[0, l]$, $\varphi(x), ~ \nu(x)\in L_2(0, l)$ and $f(x, t)\in C^{0, 1}[0, l]\times [-T, 0]$,  then the series of $u(x, t)$ is bounded with convergent numerical series with respect to $n$ and from Weierstrass M-test the series of $u(x, t)$ converges uniformly in $\Omega_2$.

 Next, we show the uniform convergence of the series representation  of $u_{xx}(x, t)$, which is given by in $\Omega_1$
 \begin{equation*}
 	u_{xx}(x, t)=-\sum\limits_{n=1}^{\infty}\lambda_n\tau_n t^{\gamma_1-1}E_{\delta_1, \gamma_1}(-\lambda_nt^{\delta_1})-
 	\end{equation*}
 	\begin{equation*}
 	-\sum\limits_{n=1}^{\infty}\lambda_n\int\limits_{0}^{t}(t-s)^{\delta_1-1}E_{\delta_1, \delta_1}\left[-\lambda_n(t-s)^{\delta_1}\right]f_n(s)ds\sin(\sqrt{\lambda_n}x),
 \end{equation*}
 and as such given in $\Omega_2$
 \begin{equation*}
 	u_{xx}(x, t)=-\sum\limits_{n=1}^{\infty}\lambda_n\varphi_n(-t)^{\gamma_2-2}E_{\delta_2, \gamma_2-1}[-\lambda_n(-t)^{\delta_2}]\sin(\sqrt{\lambda_n}x)+
 	\end{equation*}
 	\begin{equation*}
 	+\sum\limits_{n=1}^{\infty}\lambda_n\nu_n(-t)^{\gamma_2-1}E_{\delta_2, \gamma_2}[-\lambda_n(-t)^{\delta_2}]\sin(\sqrt{\lambda_n}x)-
 \end{equation*}
 \begin{equation*}
 	-\sum\limits_{n=1}^{\infty}\lambda_n\int\limits_{t}^{0}(z-t)^{\delta_2-1}E_{\delta_2, \delta_2}[-\lambda_n(z-t)^{\delta_2}]f_n(z)dz\sin(\sqrt{\lambda_n}x).
 \end{equation*}

 Integrating by parts on second part and considering Lemma \ref{lm1} we have the following estimates $~ t>0$
 $$
 |u_{xx}(x, t)|\leq
 $$
 $$
 \leq\sum\limits_{n=1}^{\infty}\left(\frac{|\tau_{2n}||t|^{\gamma_1-1}|M}{1+\lambda_n|t^{\delta_1}|}
 +\frac{|f_{2n}(0)||t|^{\delta_1}}{1+\lambda_n|t|^{\delta_1}}+\int\limits_{0}^{t}\frac{|t-s|^{\delta_1}M}{1+\lambda_n|t-s|^{\delta_1}}|f'_{2n}(s)|ds\right),
 $$

 and as such for $~ t<0$
 \begin{equation*}
 	|u_{xx}(x, t)|\leq \sum\limits_{n=1}^{\infty}\left(\frac{\lambda_n\varphi_n|-t|^{\gamma_2-2}}{1+\lambda_n|-t|^{\delta_2}}+\frac{\lambda_n\nu_n|-t|^{\gamma_2-1}}{1+\lambda_n|-t|^{\delta_2}}\right)+
 \end{equation*}
 $$
 +\sum\limits_{n=1}^{\infty}\left(\frac{|f_{2n}(0)||-t|^{\delta_2} M}{1+\lambda_n |-t|^{\delta_2}}+\int\limits_{t}^{0}\frac{|z-t|^{\delta_2} M}{1+\lambda_n |z-t|^{\delta_2}} |f'_{2n}(s)|ds\right),
 $$
 where
 $\tau_n=\frac{\tau_{2n}}{\lambda_n}$, $\varphi_n=\frac{\varphi_{2n}}{\lambda_n}$, $\nu_n=\frac{\nu_{2n}}{\lambda_n}$,

 \begin{equation*}\label{fn}
 	f_{2n}(t)=\frac{2}{l}\left\{\begin{array}{l} \int\limits_{0}^{t}f_{xx}(x, t)\sin(\sqrt{\lambda_n}x),  \,\,t>0,\\ \int\limits_{t}^{0}f_{xx}(x, t)\sin(\sqrt{\lambda_n}x), \,\,t<0.\end{array}\right.
 \end{equation*}

 If  $f(x, t)\in C^{2, 1}(0, l)\times(-T, T)$ and $\tau(x), \varphi(x), \nu(x)\in C^2(0, l)$ and  $\tau'''(x), \varphi'''(x), \nu'''(x)\in L_2(0, l)$ which are required for $t\to 0$, then, the series representation of $u_{xx}(x, t)$ is bounded with the convergent numerical series and from Weierstrass M-test the series of $u_{xx}(x, t)$ converges uniformly in $\Omega_1\cup\Omega_2$.

  Finally, the uniform convergence of the series representation of $D_{0\pm}^{(\alpha_i, \beta_i)\mu_i}u(x, t)$, $i=\overline{1,2}$ can be done similarly to the convergence of the series of  $u_{xx}(x, t)$ considering Eq.(\ref{eq1}).

 Moreover, according to (\ref{tau})-(\ref{fi}) we can see that $\tau(x), \varphi(x)$ and $ \nu(x)$  functions are written in terms of the given functions $\psi(x)$ and $f(x, t)$. For that reason we write sufficent conditions for those given functions   in order to show that all imposed conditions for $\tau(x)$, $\varphi(x)$ and $\nu(x)$ are valid, i.e

 $$\tau(x), \varphi(x), \nu(x)\in C[0, l], ~  ~ \tau(x), \varphi(x), \nu(x)\in C^2(0, l)~ ~ \text{and} ~ ~$$

 $$ \tau'''(x), \varphi'''(x), \nu'''(x)\in L_2(0, l), ~ ~ f(x, t)\in C[0, l]\times[-T, T],$$
 $$ f(x, t)\in C^{2, 1}(0, l)\times(-T, T).
 $$
If we find sufficent conditions for given functions in order to show the validity conditions of $\nu(x)$, it can be clearly seen that those sufficient conditions can be considered enough for showing that conditions for $\tau(x), \varphi(x)$ are also valid automatically. Hence we have the following equality from (\ref{nu})

 \begin{equation*}
 	\nu_n=\frac{-\lambda_n}{\Delta_n}(\psi_n+F_n)=-\frac{1}{\Delta_n\lambda_n\sqrt{\lambda_n}}\psi_{5n}-\frac{1}{\Delta_n{\lambda_n\sqrt{\lambda_n}}}F_{3n},	
 \end{equation*}
 Since the given functions can be written in the form of a Fourier series and the last equality we have the following conditions for the given functions
 $$
 \psi(x)\in C[0, l]\cap C^{4}(0, l) ~ \text{and} ~ \psi^{(5)}(x)\in L_2(0, l),
 $$
 $$
 f(x, t)\in C[0, l]\times[-T, T]\cap C^{2, 1}(0, l)\times(-T, T) ~ \text{and} ~ f_{x}^{(3)}(\cdot, t)\in L_2(0, l),
 $$
 where we assume that $\Delta_n\neq0$, $\psi(0)=\psi(l)=0$, $\psi''(0)=\psi''(l)=0$, $\psi^{(4)}(0)=0$, \\$\psi^{(4)}(l)=0$, ~ $f(0,t)=f(l, t)=f_{xx}(0, t)=f_{xx}(l, t)=0$ and we have used the following inequality
 $$2|\frac{1}{\Delta_n\sqrt{\lambda_n}}\psi_{5n}|\leq\frac{1}{\Delta_n^2 \lambda_n}+|\psi_{5n}|^2,$$
and Parseval's identity
 $$
 \sum_{n=1}^{\infty}|\psi_{5n}|^2=\|\psi^{(5)}(x)\|^2,
 $$

 \begin{equation*}
 	\psi_{n}^{(5)}=\frac{2}{l}\int\limits_{0}^{l}\psi^{(5)}(x)\sin(\sqrt{\lambda_n}x)dx,	
 \end{equation*}
 \begin{equation*}
 	F_{3n}=	\int\limits_{0}^{T}(T-z)^{\delta_1-1}E_{\delta_1, \delta_1}(-\lambda_n(T-z)^{\delta_1})f_{3n}(z)dz-
 	\end{equation*}
 	\begin{equation*}
 	-\int\limits_{-T}^{0}(z+T)^{\delta_2-1}E_{\delta_2, \delta_2}(-\lambda_n(z+T)^{\delta_2})f_{3n}(z)dz,
 \end{equation*}
$$
|F_{3n}(t)|\leq|f_{3n}(0+)|\frac{MT^{\delta_1}}{1+\lambda_nT^{\delta_1}}+\int\limits_{0}^{T}|T-z|^{\delta_1}\frac{M}{1+\lambda_n|T-z|^{\delta_1}}|f'_{3n}(z)|dz+
$$
$$
+|f_{3n}(0-)|T^{\delta_2}\frac{M}{1+\lambda_nT^{\delta_2}}+\int\limits_{-T}^{0}|z+T|^{\delta_2}\frac{M}{1+\lambda_nT^{\delta_2}}|f'_{3n}(z)|dz
$$
 \begin{equation*}
 	f_{3n}(t)=\frac{2}{l}\int\limits_{0}^{l}f^{(3)}_{x}(x, t)\cos(\sqrt{\lambda_n}x)dx.
 \end{equation*}
\begin{equation*}
	f_{3n}(0)=\frac{2}{l}\int\limits_{0}^{l}f^{(3)}_{x}(x, 0)\cos(\sqrt{\lambda_n}x)dx.
\end{equation*}

 All in all, we have just proved the following theorem (\ref{thm1}).

\begin{thm}\label{thm1} Assume that the following conditions hold:
	
	$\Delta_n\neq0$;
	
	$\psi(x)\in C[0, l]\cap C^{4}(0, l) $ such that $\psi(0)=\psi(l)=0$, $\psi''(0)=\psi''(l)=0$, $\psi^{(4)}(0)=\psi^{(4)}(l)=0$  and  $\psi^{(5)}(x)\in L_2(0, l)$; 
	
	$f(x, t)\in C[0, l]\times[-T, T]\cap C^{2, 1}(0, l)\times(-T, T)$   such that $f(0,t)=f(l, t)=0$, $f_{xx}(0, t)=f_{xx}(l, t)=0$,   $f_{x}^{(3)}(\cdot, t)\in L_2(0, l)$;
	
	then, there exists the unique solution of the considered problem.
	
\end{thm}

\bigskip

\bibliographystyle{unsrt}  


\end{document}